\documentclass[graybox]{svmult}
\usepackage{newtxtext,newtxmath}
\usepackage{siunitx}
\usepackage{helvet}         % selects Helvetica as sans-serif font
\usepackage{courier}        % selects Courier as typewriter font
\usepackage{makeidx}         % allows index generation
\usepackage{graphicx}        % standard LaTeX graphics tool
\usepackage[bottom]{footmisc}% places footnotes at page bottom
\usepackage{url}
\makeindex            

\usepackage{wrapfig,lipsum}
\begin{document}

\title*{\textbf{ Topological surgery in the small and in the large}} 
\author{Stathis Antoniou, Louis H.Kauffman and Sofia Lambropoulou}
\institute{Stathis Antoniou \at {School of Applied Mathematical and Physical Sciences \\ National Technical University of Athens, Greece} \\ \email{santoniou@math.ntua.gr}
\and Louis H.Kauffman  \at {Dept of Mathematics, Statistics, and Computer Science, University of Illinois at Chicago, USA } \\ Dept of Mechanics and Mathematics, Novosibirsk State University, Russia \\  \email{kauffman@uic.edu} 
\and Sofia Lambropoulou  \at School of Applied Mathematical and Physical Sciences \\ National Technical University of Athens, Greece \\ \email{sofia@math.ntua.gr}}
\maketitle

%%%%%%%%%%%%%%%%%%%%%%%%%%%%%%%%%%%%%%%%%%%%%%%%%%%%%%%%%%%%%%%%%%%%%%%%%%%%%%%%%%
\abstract{We directly connect topological changes that can occur in mathematical three-space via surgery, with black hole formation, the formation of wormholes and new generalizations of these phenomena. This work widens the bridge between topology and natural sciences and creates a new platform for exploring geometrical physics.}

\let\thefootnote\relax\footnotetext{

{\noindent}\textit{2010 Mathematics Subject Classification}: 57M25, 57R65, 83F05 \\
{\noindent}\textit{Keywords}: topological surgery, topological process, three-space, three-sphere, three-manifold, handle, Poincar\'{e} dodecahedral space, knot theory, natural phenomena, natural processes, dynamics, reconnection, Morse theory, mathematical model, Falaco solitons, black holes, wormholes, Einstein-Rosen bridge, cosmic string, quantum gravity, cosmology, ER=EPR, entanglement, DNA recombination, biology}
 
%%%%%%%%%%%%%%%%%%%%%%%%%%%%%%%%%%%%%%%%%%%%%%%%%%%%%%%%%%%%%%%%%%%%%%%%%%%%%%%%%%
%\newpage
\section{Introduction}
The universe undergoes topological and geometrical changes at all scales. This paper  goes to the foundations of these changes by offering a novel topological perspective.  The common features of these changes are described via topological surgery, a manifold-changing process which has been used in the study and classification of manifolds. We briefly address small and large scale phenomena exhibiting $1$ and $2$-dimensional surgery and then focus on large scale cosmic phenomena exhibiting $3$-dimensional surgery. More precisely, we describe the formation of black holes and wormholes. Our surgery approach allows the formation of a black hole from the collapse of a knotted cosmic string, without ending up in a singular manifold. It further describes Einstein-Rosen bridges (wormholes) linking the two black holes through a singularity where the disconnected black holes collapse to each other, and the bridge is born topologically. The collapse of a cosmic string can be viewed as an orchestrated creation of bridges that is topologically equivalent to $3$-dimensional surgery. We present a rich family of $3$-manifolds that can occur and the possible implications of these constructions in quantum gravity and general relativity.

\section{The topological process of surgery}
Topological surgery is a mathematical technique introduced by A.H.Wallace~\cite{Wal} and J.W.Milnor~\cite{Milsur} which creates new manifolds out of known ones in a controlled way. Given an $m$-manifold $M$, an $m$-dimensional $n$-surgery consists of removing a thickened sphere $S^{n} \times D^{m-n}$ and gluing back another thickened sphere $D^{n+1} \times S^{m-n-1}$ using a gluing homeomorphism along the common boundary $S^{n} \times S^{m-n-1}$, see~\cite{Ra} for details. This operation produces a new $m$-manifold $M'$ which may, or may not, be homeomorphic to $M$. Since $(S^{n} \times D^{m-n}) \cup (D^{n+1} \times S^{m-n-1})=\partial (D^{n+1}\times D^{m-n}) \cong D^{m+1}$, an $m$-dimensional $n$-surgery can be seen as the process of passing from one boundary component of handle $D^{m+1}$ to the other. The extra dimension of the $(m+1)$-dimensional handle leaves room for continuously passing from one boundary component of the handle to the other. 

\begin{figure}[!h]
%\sidecaption
%\includegraphics[scale=.1366]{SurgEx.jpg}
\includegraphics[scale=.19]{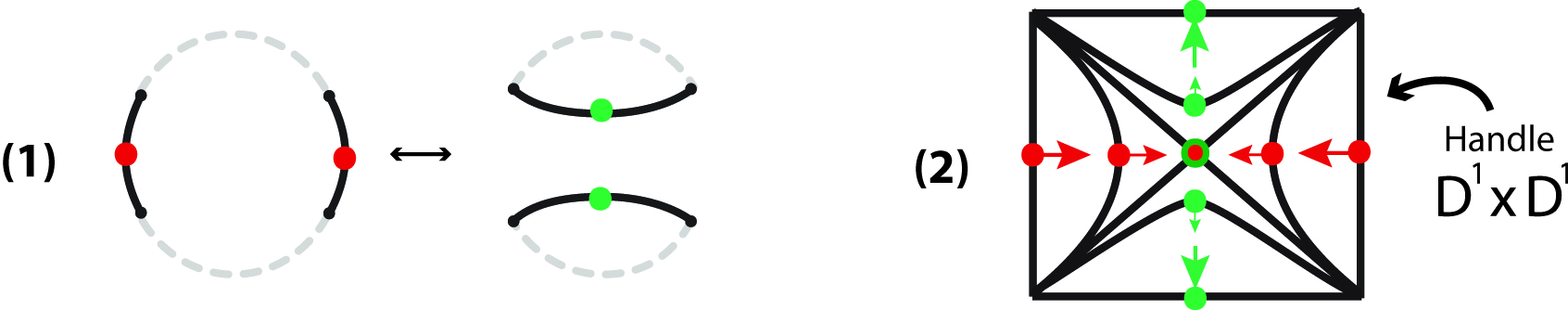}
\caption{\textbf{(1)} $1$-dimensional surgery \textbf{(2)} $2$-dimensional handle }
\label{SurgEx}
\end{figure}

For example, the process of $1$-dimensional $0$-surgery shown in Fig.~\ref{SurgEx}~(1) removes the $1$-dimensional thickening of a $0$-sphere (represented by the two red points) and replaces it with the thickening of another $0$-sphere (represented by two green points). Starting with the circle $M=S^1$, this process produces two circles $M'=S^0 \times S^1$. This global change of topology is induced by the local process of collapsing the cores of two segments (the two red points) and uncollapsing the core of the other two segments (the two green points), see Fig.~\ref{SurgEx}~(2). As also shown in the figure, this local process happens within a $2$-dimensional handle $D^1 \times D^1$. One can provide a algebraic description of this dynamic process by using the local form of a Morse function, see Lemma 2.19 of~\cite{Ra}. For instance, in dimension $1$, the local process is described by varying parameter $t$ of the level curves $-x^2+y^2=t$ in the range $t \in (-1,1)$.

\section{Small scale surgery in Nature}
We will briefly present some natural phenomena exhibiting topological surgery in dimensions $1$ and $2$. All these phenomena happen in small scales, meaning that their characteristic lengths range from the size of a molecule to a few meters. 
 
\begin{figure}[!h]
\sidecaption
\includegraphics[scale=.12]{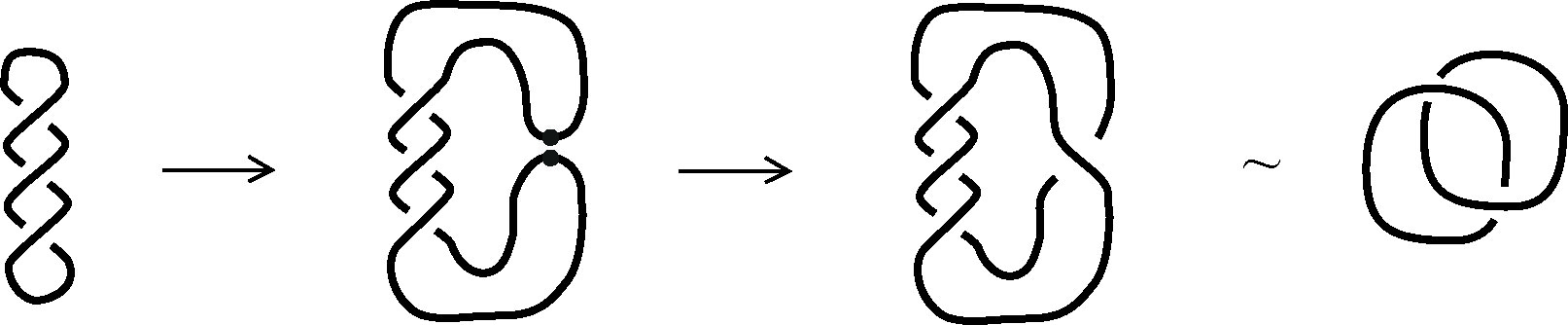}
\caption{DNA recombination as an example of $1$-dimensional $0$-surgery}
\label{DNArecomb}
\end{figure}

The process of $1$-dimensional surgery is exhibited in various natural phenomena where segments are detached and rejoined such as the crossing over of chromosomes during meiosis, viscous vortex reconnection and site-specific DNA recombination. These phenomena have been detailed in~\cite{SS1,SS2,SS3} where we show that although they are quite different, they undergo a similar topological change which is described using our surgery approach.

If the initial manifold is an embedding of the circle, $1$-dimensional $0$-surgery can create or destroy a crossing hence producing new knots or links. For example, starting with the circular DNA molecule of Fig.~\ref{DNArecomb}, with the help of certain enzymes, site-specific recombination performs a $1$-dimensional $0$-surgery on the molecule and produces the Hopf link. 

\begin{figure}[!h]
\sidecaption
\includegraphics[scale=.065]{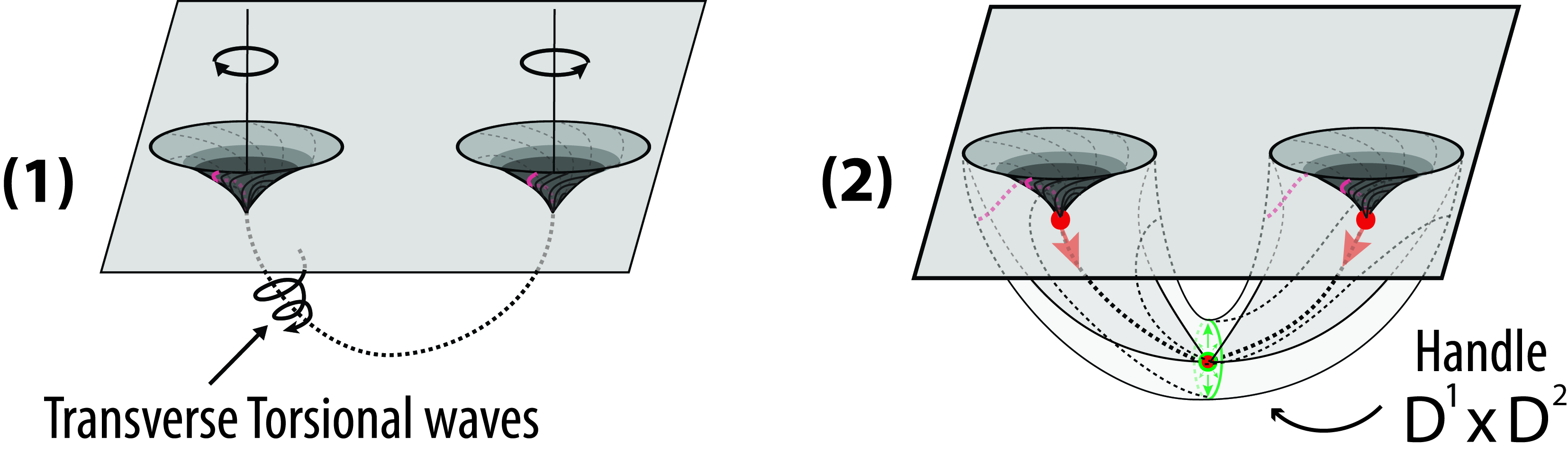}
\caption{\textbf{(1)} Falaco topological defects \textbf{(2)} 3-dimensional handle}
\label{2D_Falacohomeo2}
\end{figure}

Nature is filled with $2$-dimensional surgeries too, see~\cite{SS1,SS2,SS3} for details. Examples comprise gene transfer  in bacteria, where the donor cell produces a connecting tube called a `pilus' which attaches to the recipient cell, the biological process of mitosis,  where a cell splits into two new cells, and the formation Falaco Solitons. We will describe how surgery is exhibited in Falaco Solitons, the dynamics of which are visible to the naked eye. Each Falaco soliton consists of a pair of contra-rotating identations in the water-air surface of a swimming pool, see Fig.~\ref{2D_Falacohomeo2}~(1) and \cite{Ki}. From the topological viewpoint the surgery consists in taking disk neighborhoods of two points $S^0 \times D^2$ (the identations in Fig.~\ref{2D_Falacohomeo2}~(1)) and joining them via a tube (which is a thickened circle $D^1 \times S^1$), see Fig.~\ref{2D_Falacohomeo2}~(2). Here the tube is the cylindrical vortex made from the propagation of the torsional waves around the singular thread. The $3$-dimensional handle containing all the $2$-dimensional temporal `slices' of this process is shown Fig.~\ref{2D_Falacohomeo2}~(2).

\section{Large scale surgery in Nature}
Large scale phenomena can also exhibit $1$- or $2$-dimensional surgery. For instance,  $1$-dimensional surgery happens in magnetic reconnection, the phenomena whereby cosmic magnetic field lines from different magnetic domains are spliced to one another, changing their pattern of conductivity with respect to the sources, see~\cite{SS1} for details. Moving up to $3$-dimensional surgery, we have two types of surgery, both of which require four dimensions to be visualized. As we will see, both types describe large scale cosmic phenomena. The first type is exhibited in the creation of entangled black holes while the second one describes the formation of black holes from cosmic strings.

\subsection{Types of 3-dimensional surgery}
Starting with a 3-manifold $M$, we first have the \textit{ $3$-dimensional $0$-surgery}, whereby two 3-balls $S^0\times D^3$ are removed from  $M$ and are replaced in the closure of the remaining manifold by a thickened sphere $D^1\times S^2$:

\begin{samepage} 
 \begin{center}
$\chi(M) = \overline{M\setminus h(S^0\times D^{3})} \cup_{h} (D^{1}\times S^{2})$
\end{center}  
\end{samepage} 

Next, for $m=3$ and $n=2$, we have the \textit{ $3$-dimensional $2$-surgery}, which is the reverse (dual) process of $3$-dimensional $0$-surgery. Hence we will not consider it a different type of $3$-dimensional surgery.

Finally, for $m=3$ and $n=1$, we have the \textit{ $3$-dimensional $1$-surgery}, whereby a solid torus $S^1\times D^2$ is removed from $M$ and is replaced by another solid torus $D^2\times S^1$ (with the factors now reversed) via a homeomorphism $h$ of the common boundary:

\begin{samepage} 
 \begin{center}
$\chi(M) = \overline{M\setminus h(S^1\times D^{2})} \cup_{h} (D^{2}\times S^{1})$
\end{center}  
\end{samepage} 

This type of surgery is clearly self-dual.

\subsection{3-dimensional 0-surgery and entangled black holes}
The process of $3$-dimensional $0$-surgery joins the spherical neighborhoods of two points via a tube $D^1 \times S^2$ which is one dimension higher than the one shown in Fig.~\ref{2D_Falacohomeo2}~(2). If we consider that our initial manifold is the $3$-dimensional spatial section of the $4$-dimensional spacetime, this tube is what physicists call a wormhole. A connection between Falaco solitons and wormholes has been conjectured by R.M. Kiehn~\cite{KiehnSmall}. Our surgery description reinforces this connection. Moreover, this change of topology, which, according to J.A. Wheeler, results from quantum fluctuations at the Planck scale~\cite{Whee_QF}, can now also be viewed as a result of a `classical' continuous topological change of $3$-space.

\begin{figure}[!h]
\sidecaption
\includegraphics[scale=.065]{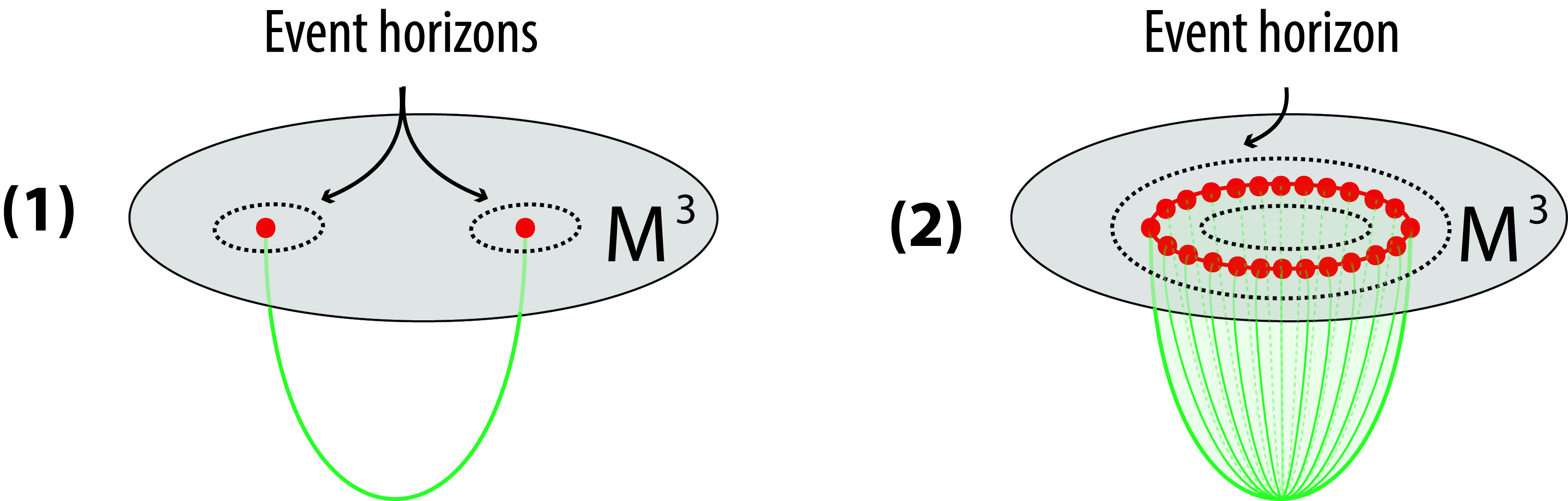}
\caption{\textbf{(1)} Pair of entangled black holes \textbf{(2)} String of entangled black holes}
\label{EnStrings}
\end{figure}

Let us now consider the $ER = EPR$ hypothesis of L. Susskind~\cite{ER_EPR}, which says that a wormhole is equivalent to the quantum entanglement of two concentrated masses that each forms its own black hole. Adding this hypothesis to our description, the two sites in space are the singularities of the two black holes, shown in red in Fig.~\ref{EnStrings}~(1), which will not collapse individually but will become the ends of the wormhole, shown in green in Fig.~\ref{EnStrings}~(1). We cannot visualize this process directly but it can be understood by considering that the green arc is the core $D^1$ of the higher dimensional handle $D^1 \times D^3$ for the wormhole. Note that an observer in our initial $3$-space $M^3$ would not be able to detect the topological change,  which occurs across the event horizons.

\subsection{3-dimensional 1-surgery and cosmic string black holes}
\begin{figure}[!h]
\sidecaption
\includegraphics[scale=0.099]{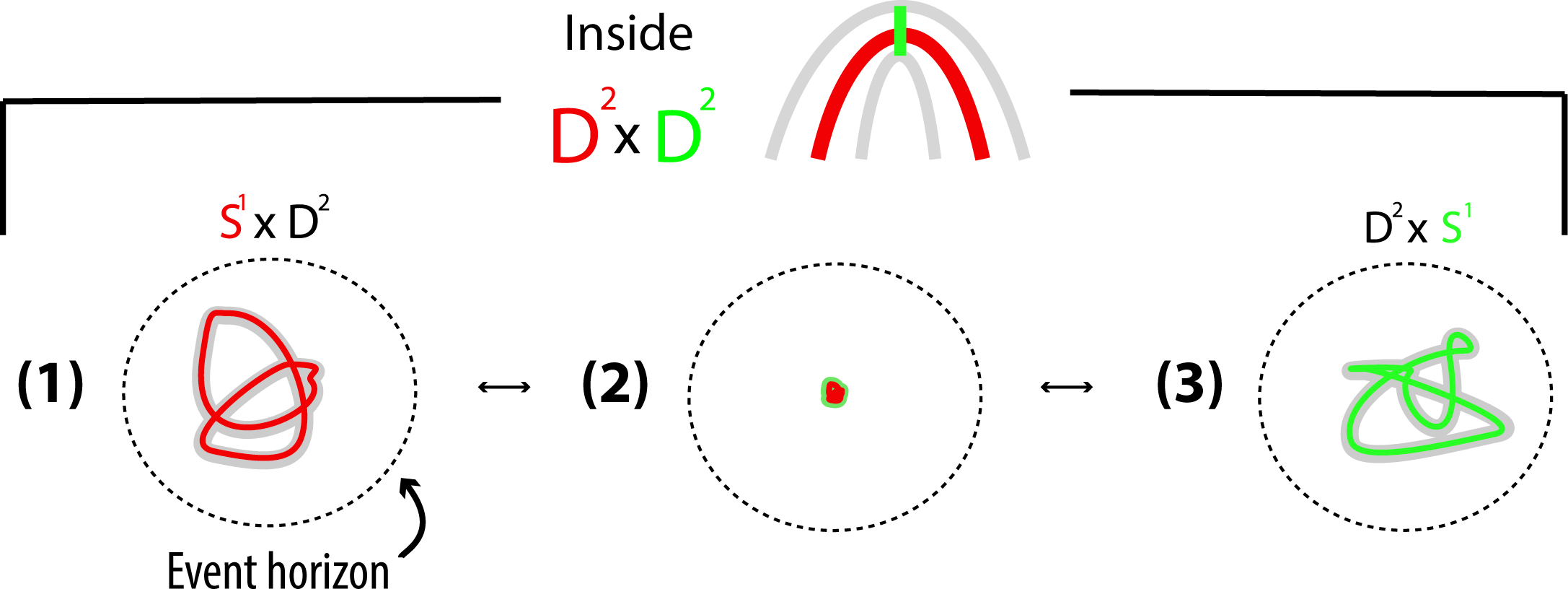}
\caption{3-dimensional 1-surgery inside the event horizon}
\label{Fig35}
\end{figure}

The other type of $3$-dimensional surgery describes a more subtle topological change. It collapses a solid torus (which is a thickened circle) to a point and uncollapses another solid torus in such way that the meridians are glued to the longitudes and vice-versa. This type of surgery is also called `knot surgery' as the circle can be a knot. Knot surgery is an ideal candidate for describing black holes that are formed via the collapse of cosmic strings. This idea is based on~\cite{Hawk} where S.W. Hawking estimates that a fraction of cosmic string loops can collapse to a small size inside their Schwarzschild radius. As cosmic strings are hypothetical topological defects of small (but non-zero) diameter, a cosmic string loop can be considered as a knotted solid torus. As described in~\cite{Hawk}, the loop collapses to a point thus creating a black hole the center of which contains the singularity. At that point, the $3$-space becomes singular, see the passage from Fig.~\ref{Fig35}~(1) to~(2).    

Our surgery description says more~\cite{AKL}. According to it the process doesn't stop at the  singularity, but continues with the uncollapsing of another cosmic string loop from the singularity, see Fig.~\ref{Fig35}~(3). Thus, the creation of a cosmic string black hole is a $3$-dimensional $1$-surgery that continuously changes the initial $3$-space to another $3$-manifold. The process goes through the singular point of the black hole without having a singular manifold in the end. Instead, one ends up with a topologically new universe with a local topology change in the $3$-space, which happens within the event horizon.

This type of surgery is also related to the $ER = EPR$ hypothesis. Consider a cosmic string made of pairs of entangled concentrated masses. When each pair of masses collapse, they become connected by a wormhole as previously shown in Fig.~\ref{EnStrings}~(1). Given that all these pairs of masses have started on the same cosmic string, the distinct wormholes merge and the entire collection of wormhole cores (the green arcs in Fig.~\ref{EnStrings}~(1)) forms a 2-disc $D^2$, see Fig.~\ref{EnStrings}~(2), which is the core of the higher dimensional handle $D^2 \times D^2$ containing the temporal `slices' of the process. Our surgery description generalizes having a separate Einstein-Rosen bridge for each pair of black holes and amalgamates these bridges to form a new manifold in three dimensions. The effect of surgery is that, from any black hole location on the cosmic string to any other, there is a `bridge' through the new $3$-manifold. As this process joins the neighborhood of a circle instead of two points, one can rotate Fig.~\ref{EnStrings}~(1) to receive Fig.~\ref{EnStrings}~(2).

\begin{figure}[!h]
\sidecaption
\includegraphics[scale=.16]{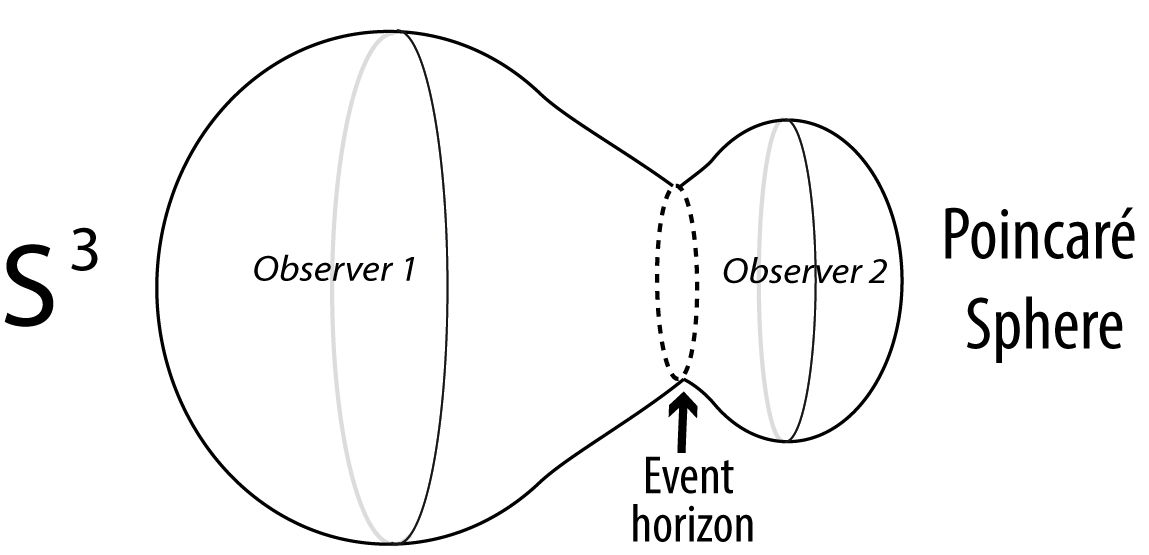}
\caption{Observer 1 in $S^3$ and  Observer 2 in the Poincar\'{e} dodecahedral space}
\label{BH_PoincaFlip}
\end{figure}

Another advantage of knot surgery is that it is able produce a great variety of $3$-manifolds. In fact, according to a theorem by A.H. Wallace~\cite{Wal} and W. Lickorish~\cite{LickTh} knot surgery can create all closed, connected, orientable $3$-manifolds. One such $3$-manifold, which is of great interest to physicists, is the Poincar\'{e} dodecahedral space, that has been proposed as possible shape for the geometric universe~\cite{Weeks,Luminet,Levin}. This manifold can be obtained by doing knot surgery on the trefoil knot (with the right framing, see~\cite{PS}). Did the shape of the universe come about via the collapse of a trefoil cosmic string?!

Suppose there are observers in an initial spherical universe $M^3=S^3$ containing a trefoil cosmic string. After surgery, a `mathematical' observer would be able to see the Poincar\'{e} dodecahedral space and detect the topology change.  However, a physical observer, who is subject to the restrictions of physical laws, would only see towards the event horizon in which the trefoil cosmic string has collapsed. Let us call this observer, \textit{Observer 1}, see Fig.~\ref{BH_PoincaFlip}. After surgery, \textit{Observer 1} would see the same universe $S^3$, the only change being the formation of the event horizon. On the other side of the event horizon, a new universe emerges in which new observers might evolve. Such an  observer, say \textit{Observer 2}, will see a Poincar\'{e} dodecahedral space and the event horizon from the other side, unaware that the original $S^3$ universe is behind it, see Fig.~\ref{BH_PoincaFlip}. Finding the Poincar\'{e} dodecahedral space (or some other non-trivial $3$-manifold) in our universe may indicate that we are observers that evolved inside the event horizon of a collapsed trefoil cosmic string (or some other cosmic string).

\section{Conclusion}
The surgery approach provides continuous paths to wormhole and cosmic string black hole formation. If one adds the $ER = EPR$ hypothesis, surgery also describes the entanglement of a pair or a string of black holes. Our topological perspective offers a process producing black holes and new non-singular $3$-manifolds from cosmic strings, binding entanglement and the connectivity of space with the rich structure of three- and four-dimensional manifolds.
 
\begin{acknowledgement} Kauffman's work was supported by the Laboratory of Topology and Dynamics, Novosibirsk State University (contract no. 14.Y26.31.0025 with the Ministry of Education and Science of the Russian Federation).
\end{acknowledgement}

%%%%%%%%%%%%%%%%%%%%%%%%%%%%%%%%%%%%%%%%%%%%%%%%%%%%%%%%%%%%%%%%%%%%%%%%%%%%%%%%%
%\newpage

\end{document}